\documentclass[12pt,a4paper]{article}
\usepackage[12pt]{extsizes}
\usepackage[T2A]{fontenc}
\usepackage[cp1251]{inputenc}
\usepackage[english]{babel}
\usepackage{amssymb,amsfonts,amsmath,mathtext}
\usepackage{cite,enumerate,float,indentfirst}
\usepackage{tikz}
\usepackage{calc}
\usepackage{fancyhdr}
\usepackage{textcase}
\usepackage{amsthm}
\usepackage{tocloft}
\usepackage[left=3cm,right=2cm,top=2cm,bottom=2cm,bindingoffset=0cm]{geometry}

\makeatletter
\renewcommand{\@biblabel}[1]{#1.}
\makeatother

\newcommand{\Log}{\mathop{\mathrm{Log}}\nolimits}

\title{A List of Integral Representations for Diagonals\\ of Power Series of Rational Functions}
\author{Artem Senashov\footnote{The author was supported by the Foundation for the Advancement of Theoretical Physics and Mathematics "BASIS"(18-1-7-60-3).}}

\date{}

\begin{document}

	\maketitle
	\begin{abstract}

 In this article, we present the integral representations of the power series diagonals.
Such representations are obtained by lowering the integration multiplicity for the previously known integral representation. The procedure is carried out within the framework of Leray's residue theory. The concept of the amoeba of the complex analytical hypersurface plays an essential role in the construction of new integral representations.
\end{abstract}

\bigskip

\section{Introduction}


 A range of problems associated with branching of parametric integrals is concerned with a study of the diagonals of power series \cite{Pochecutov dis}, \cite{DenefLipshitz} and~\cite{Lipshitz}.
It should be noted that much earlier the concept of the diagonal of a power series was used by A. Poincare \cite{Puankare} to study the anomalies of planetary motion.

{\it The diagonal} of a Laurent power series
\begin{equation}\label{b1}
F(z)=\sum_{\alpha\in \mathbb{Z}^n}c_{\alpha}z^{\alpha}
\end{equation}
is defined to be the generating function of a subsequence of coefficients $\{c_{\alpha}\}_{\alpha\in L}$ indexed by elements $\alpha$ of some sublattice  $L\subset \mathbb{Z}^n$ (see   \cite{Pochecutov dis} and \cite{Pochecutov}). Such diagonals are called complete. Diagonals are graded according to the dimension (rank) of the sublattice.

Following \cite{Pochecutov dis}, we describe the specifics of the problem on the properties of the diagonals of series for rational functions of $ n $ variables

\begin{equation}\label{b3}
F(z)=\frac{P(z)}{Q(z)}=\frac{P(z_1,\dots,z_n)}{Q(z_1,\dots,z_n)},
\end{equation}
where $P$ and $Q$ are irreducible polynomials. Consider an arbitrary Laurent series for $ F $ centered at zero:

$$
F(z)=\sum_{\alpha\in \mathbb{Z}^n}c_{\alpha}z^{\alpha}=\sum_{\alpha\in \mathbb{Z}^n}c_{\alpha_1,\dots,\alpha_n}z_1^{\alpha_1}\dots z_n^{\alpha_n}.
$$

It is known that such a series converges in the domain $\Log^{-1}(E)$, where $E$ is a connected component of the complement $R^n\backslash A_Q $ of amoeba of the denominator $Q$ \cite{GKZ}. Recall that the {\it amoeba} $A_Q$ of the polynomial $Q$ or of the algebraic hypersurface
$$
V=\{z\in (\mathbb{C}\backslash0)^n:Q(z)=0\}
$$
is defined to be the image of $ V $ under the mapping $\Log:(\mathbb{C}\backslash0)^n\rightarrow \mathbb{R}^n$, given as follows
 $$
 \Log:(z_1,\dots,z_n)\rightarrow(\log|z_1|,\dots,\log|z_n|).
 $$
According to a result of the article \cite{FPT}, there is an injective order function

 $$\nu:{E}\rightarrow \mathbb{Z}^n\bigcap N_Q,$$
mapping each connected component $ E $ of the complement $\mathbb {R}^n\backslash A_Q$ to integer vector $\nu=\nu(E)$, belonging to the Newton polytope $N_Q$ of the polynomial $Q$. Thus, all connected components can be indexed as $\{E_\nu\}$, where $\nu$ runs over some subset of integer points in $N_Q$. For example, the Taylor series of a rational function $\frac{P}{Q},Q(0)\neq0$ converges in the preimage $\Log^{-1}(E_0)$ of the component $E_0$.

Let us consider in more detail the $p$-dimensional diagonal of the series (\ref{b1}). For this, in $\mathbb{Z}^n_+$ we fix a basis $q^{(1)},\dots, q^{(p)}$ of the sublattice $l\subset L$. Remark, that these vectors define a basis of $l$ if and only if they are relatively prime. It is equivalent to the property that the set $q^{(1)},\dots, q^{(p)}$ can be complemented to a unimodular matrix by integer vectors $q^{(p+1)},\dots,q^{(n)}$ (see \cite{NPT} or \cite[Proposition 4.2.13]{SadTsikh}). Denote this matrix by
$$
A=(q^{(1)},\dots,q^{(n)})
$$
and recall that its determinant is equal to 1. Directions $q^{(1)},\dots, q^{(p)}$ define a diagonal subsequence $\{c_{lq}\}_{l\in \mathbb{Z}^n_+}$ of the  multiple sequence of coefficients $\{c_{\alpha}\}$ of the series (\ref{b1}).

The generating function
$$
d_q(t)=d(t_1,\dots,t_p)_{q_1,\dots,q_p}=\sum_{l \in \mathbb{Z}^p_+}c_{l_1q^{(1)}+\dots+l_pq^{(p)}}t_1^{l_1}\cdots t_p^{l_p}$$
of the subsequence $\{c_{lq}\}_{l\in \mathbb{Z}^n_+}$ is called {\it the one-sided $ q $-diagonal} of the series (\ref{b1}).

If $t\in \Log^{-1}(E_0)$, the $p$-dimensional diagonal $d_q(t)$ of the Taylor series admits the integral representation

\begin{equation}\label{c1}
d_q(t)=\frac{1}{(2\pi i)^n}\int_{\Gamma_{\rho}}F(z)\frac{z^{q^{(1)}}\cdots z^{q^{(p)}}}{(z^{q^{(1)}}-t_1)\dots(z^{q^{(p)}}-t_p)}\frac{dz_1}{z_1}\dots\frac{dz_n}{z_n},
\end{equation}
where $z^q$ is a monomial $z_1^{q_1}\dots z_n^{q_n}$, and cycle

$$
\Gamma_{\rho}=\{z\in \mathbb{C}^n:|z_1|=e^{{\rho}_1},\dots,|z_n|=e^{{\rho}_n} \}
$$
is chosen so that

a) the closed polydisc

$$
\overline{U_{\rho}}=\{z\in \mathbb{C}^n:|z_1|\leq e^{{\rho}_1},\dots,|z_n|\leq e^{{\rho}_n} \}
$$

 doesn't contain poles of $F(z)$;

b) parameters $t=(t_1,\dots,t_p)$ satisfy the inequalities  $|t_i|<e^{\langle q_i,\rho\rangle},i=1,\dots,p$.

The cycle $\Gamma_{\rho}$ is a preimage $Log^{-1}\rho$ of the point $\rho$ from the connected component $E_0$ of the amoeba complement. Here we prove that the integral which represents the diagonal $ d_q (t) $ admits a decrease of the multiplisity of integration while preserving the rationality of the integrand.

 We will assume that $N_Q\subset \mathbb{R}^n_u$, and the image $A^{-1}(N_Q)\subset \mathbb{R}^n_v$, here $\mathbb{R}^n_v$ and $\mathbb{R}^n_u$ are the $n$-dimensional real variable spaces $u$ and $v$ respectively. Let us denote by $N'$  the projection of the polyhedron $A^{-1}N_Q$ on the coordinate $(n-p)$-dimensional plane $\{v\in \mathbb{R}^n: v_1=0,\dots,v_p=0\}$, and by $Q'(t,w')$ the Laurent polynomial $Q[(t,w')^{A^{-1}}] $ in variables $w'=(w_{p+1},\dots,w_n)$, wherein $t_1,\dots,t_p$ are parameters.

\noindent\textbf{Theorem 1.} \emph{The diagonal $d_q(t)$ given by (\ref{c1}) is representable by an integral in the $(n-p)$-dimensional complex algebraic torus $(\mathbb{C} \backslash 0)^{n-p}$ in variables $w_{p+1},\dots,w_{n}$ according to the formula
\begin{equation}
 d_q(t)=\frac{1}{(2\pi i)^{n-p}}\int_{\Log^{-1}({\rho'})}F[(t_1,\dots,t_p,w_{p+1},\dots,w_n)^{A^{-1}}]\frac{dw_{p+1}\dots dw_n}{w_{p+1}\dots w_n},
\end{equation}
where
 $$
\rho'=((A\rho)_{p+1},\dots,(A\rho)_n)
$$
   belongs to the connected component $E_0^{'}$ of the complement for the amoeba  $A_{Q^{'}}$ of the hypersurface $V'=\{w'\in (\mathbb{C} \backslash 0)^{n-p}:Q'(t,w')=0\}$.}

\section{Proof of Theorem 1}

It is assumed that the diagonal (\ref{c1}) is considered for the Taylor series of the rational function $F=\frac{P}{Q}$. Therefore, it implies that $Q(0)\neq0$, and it means that the origin is a vertex of the Newton polytope $N_Q$.

Since the determinant of the integer matrix $A$ is equal to one, the inverse matrix
\begin{equation*}
A^{-1}=
\begin{pmatrix}
b_1^{(1)}&\dots&b_n^{(1)}\\
\vdots&\ddots&\vdots\\
b_1^{(n)}&\dots&b_n^{(n)}
\end{pmatrix}=(b_j^{(i)})
\end{equation*}
is integer and its entries $b_j^{(i)}$ are algebraic complements to elements $q_j^{(i)}$.  Rows and columns of this matrix we denote by $b^{(i)}$ and $b_j$ respectively.
Let us make in (\ref{c1}) the change of variables

\begin{equation*}
z=w^{A^{-1}}=(w^{b_1},\dots,w^{b_n}),
\end{equation*}
which in detail looks as follows
\begin{equation*}
(z_1,\dots,z_n)=(w_1^{b_1^{(1)}}\cdots w_n^{b_1^{(n)}},w_1^{b_2^{(1)}}\cdots w_n^{b_2^{(n)}},\dots,w_1^{b_n^{(1)}}\cdots w_n^{b_n^{(n)}}).
\end{equation*}

First, note that $z^{q^{i}}$ passes to $w_i$

\begin{equation*}
\begin{split}
z^{q^{(i)}}=z_1^{ q_1^{(i)}}\dots z_n^{ q_n^{(i)}}=\Bigl(w_1^{b_1^{(1)}}\dots w_n^{b_1^{(n)}}\Bigr)^{ q_1^{(i)}}\dots\Bigl(w_1^{b_n^{(1)}}\dots w_n^{b_n^{(n)}}\Bigr)^{ q_n^{(i)}}=\\
=w_1^{\langle b^{(1)},q^{(i)}\rangle}\cdots w_n^{\langle b^{(n)},q^{(i)}\rangle}=w_i,
\end{split}
\end{equation*}
since $\langle b^{(i)},q^{(j)}\rangle=\delta_{ij}$, where $\delta_{ij}$ is the Kronecker symbol.

Applying the change of variables to the logarithmic differentials, we obtain

\begin{equation*}
\begin{split}
\frac{dz_i}{z_i}=
\frac{d(w_1^{b_i^{(1)}}\dots w_n^{b_i^{(n)}})}{w_1^{b_i^{(1)}}\dots w_n^{b_i^{(n)}}}=
\frac{\sum_{k=1}^n b_i^{(k)}w_1^{b_i^{(1)}}\dots w_k^{b_i^{(k)}-1} \dots w_n^{b_i^{(n)}}    dw_k}{w_1^{b_i^{(1)}}\dots w_n^{b_i^{(n)}}}.
\end{split}
\end{equation*}
Next we multiply the expressions for the logarithmic differentials $\frac{dz_i}{z_i}$, taking into account the properties of the external product: $dw_i\wedge dw_i=0$ and $dw_i \wedge  dw_j=-dw_j\wedge dw_i$.

As a result we obtain the expression
\begin{equation*}
\begin{split}
\frac{|A^{-1}|w_1^{\sum_{i=1}^n b_i^{(1)}-1}\dots w_n^{\sum_{i=1}^n b_i^{(n)}-1}dw_1\wedge  \dots\wedge  dw_n}{w_1^{\sum_{i=1}^n b_i^{(1)}}\dots w_n^{\sum_{i=1}^n b_i^{(n)}}}
     =\frac{dw_1 \wedge \dots \wedge dw_n}{ w_1\dots w_n},
     \end{split}
\end{equation*}
and the integral (\ref{c1}) takes the form

\begin{equation}\label{b5}
 d_q(t)=\frac{1}{(2\pi i)^{n}}\int_{\varphi_\sharp(\Gamma_{\rho})}F[(w_1,\dots,w_n)^{A^{-1}}]\frac{w_1\cdots w_p}{(w_1-t_1)\dots(w_p-t_p)} \frac{dw_1\wedge\dots \wedge dw_n}{w_1\dots w_n},
\end{equation}
where $\varphi_\sharp$ is the homomorphism induced by the mapping $\varphi:z\rightarrow w=z^A$.

The cycle $\Gamma_{\rho}$ is parameterized as follows

$$
\Log^{-1}(\rho)=\{z=e^{\rho+iA^{-1}\theta}:\theta\in A([0,2\pi)^n)\}.
$$
Hence,
$$
\varphi_\sharp(\Gamma_{\rho})=\{w=z^A:z\in\Gamma_{\rho}\}=\{w=e^{A\rho+AiA^{-1}\theta}\}=\Log^{-1}({A\rho}).
$$
In this way,
$$
\varphi_\sharp(\Gamma_{\rho})=\{ w:|w_1|=e^{(A\rho)_1},\dots,|w_n|=e^{(A\rho)_n} \},
$$
where $(A\rho)_i$ is the $i$-th component of the vector $A\rho$.

By the Cauchy formula we get

\begin{equation*}
 d_q(t)=\frac{1}{(2\pi i)^{n}}\int_{\varphi_\sharp(\Gamma_{\rho})}F[(w_1,\dots,w_n)^{A^{-1}}]\frac{dw_1}{w_1-t_1}\cdots\frac{dw_p}{w_p-t_p} \frac{dw_{p+1}\dots dw_n}{w_{p+1}\dots w_n}=
\end{equation*}

\begin{equation}\label{b7}
=\frac{1}{(2\pi i)^{n-p}}\int_{\Log^{-1}({\rho'})}F[(t_1,\dots,t_p,w_{p+1},\dots,w_n)^{A^{-1}}]\frac{dw_{p+1}\dots dw_n}{w_{p+1}\dots w_n},
\end{equation} \ \\
where
 $$
\rho'=((A\rho)_{p+1},\dots,(A\rho)_n)
$$
  belongs to the connected component $ E_0'$ of the complement of the amoeba $A_{Q'}$ for the hypersurface $V'=\{w'\in (\mathbb{C} \backslash 0)^{n-p}:Q'(t,w')=0\}$. Theorem 1 is proved.

Let us comment the reduction of the integral (\ref{c1}) to  (\ref{b7}). It is not difficult to see that the integrand in (\ref{c1}) admits representation in the form

$$
\frac{df_1}{f_1}\wedge\dots\wedge\frac{df_p}{f_p}\wedge\psi,
$$
where $\psi=\psi_p$ is a rational differential form of degree $n-p$, and $f_i=z^{q^{(i)}}-t_i$.
The system of binomial equations $f_1=0,\dots,f_p=0$ defines an $(n-p)$-dimensional complex torus $\mathbb{T}^{n-p}$ (embedded in the torus $\mathbb{T}^n=(\mathbb{C}\setminus0)^n$). In this case, the real torus $\Gamma_\rho$ is a $p$-fold tube over a real torus $\gamma\subset\mathbb{T}^{n-p}$(in the coordinates $w$, it is $Log^{-1}(A\rho')$). Thus, conditions for the application of the Leray residue formula are satisfied (see \cite{Tsikh A Yger A.},\cite{Aizenberg Yuzhakov}), so integrals (\ref{c1}) and (\ref{b7}) coincide.

\section{Example}
Consider an example of applying  Theorem 1 to find the integral representation of  the Taylor series diagonal defined by vectors $q_1=\begin{pmatrix}
1\\
1\\
1
\end{pmatrix}$ and
 $q_2=\begin{pmatrix}
1\\
2\\
2
\end{pmatrix}$ for the function $$F(z)=\frac{1}{1+z_1+z_2+z_3+z_2z_3}.$$

The Newton polytope of the denominator of the function $F$ looks as it is shown in Figure~1.

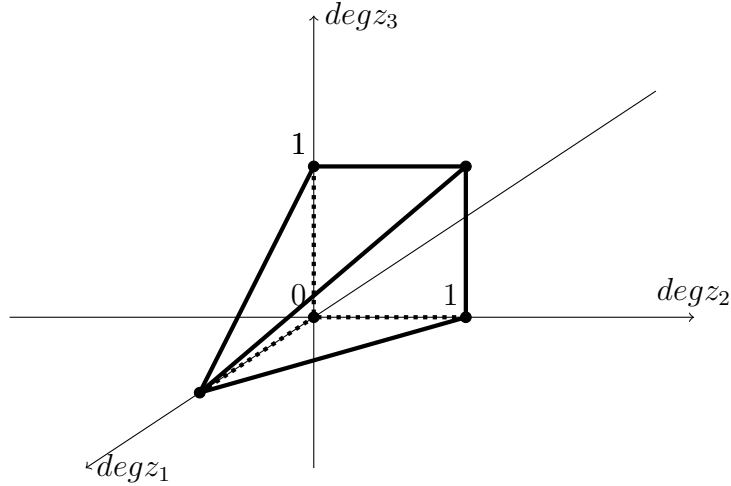
\begin{figure}[h]
\begin{center}
\begin{tikzpicture}[domain=-4:4]
\draw [->](-4,0)--(5,0) node[above]{$ deg z_2$};
\draw [->](4.5,3)--(-3,-2) node[right]{$deg  z_1$};
\draw [->](0,-2)--(0,4) node[right]{$ deg z_3$};
\draw [ultra thick](-1.5,-1)--(2,2);
\draw [ultra thick](-1.5,-1)--(2,0);
\draw [ultra thick](-1.5,-1)--(0,2) ;
\draw [ultra thick](0,2)--(2,2);
\draw [ultra thick](2,2)--(2,0);
\draw [dotted][ultra thick](-1.5,-1)--(0,0) ;
\draw [dotted][ultra thick](2,0)--(0,0) ;
\draw [dotted][ultra thick](0,2)--(0,0) ;
\draw [fill] (0,2) circle (2pt) (-1.5,-1) circle (2pt) (0,0) circle (2pt) (2,2) circle (2pt)(2,0) circle (2pt);
\draw (-0.2,0.3) node{$ 0$};
\draw (-0.2,2.3) node{$ 1$};
\draw (1.8,0.3) node{$ 1$};
\draw (-0.2,2.3) node{$ 1$};
\end{tikzpicture}
\caption{The Newton polytope of $Q(z)=1+z_1+z_2+z_3+z_2z_3$.}
\end{center}
\end{figure}

For the two-dimensional diagonal $d_{q_1,q_2}(t_1,t_2)=\sum_{l \in \mathbb{Z}^2_+}c_{l_1q^{(1)}+l_2q^{(2)}}t_1^{l_1}t_2^{l_2}$ in the set $\Log^{-1}(E_0)$, one has the following integral representation

\begin{equation}\label{b8}
d_q(t_1,t_2)=\frac{1}{(2\pi i)^3}\int_{\Gamma_{\rho}}\frac{1}{1+z_1+z_2+z_3+z_2z_3}\cdot
\frac{z_1z_2z_3\cdot z_1z_2^2z_3^2}{(z_1z_2z_3-t_1)(z_1z_2^2z_3^2-t_2)}\cdot
\frac{dz_1}{z_1}\frac{dz_2}{z_2}\frac{dz_3}{z_3},
\end{equation}
where the cycle

$$
\Gamma_{\rho}=\{z\in \mathbb{C}^n:|z_1|=e^{{\rho}_1},|z_2|=e^{{\rho}_2},|z_3|=e^{{\rho}_3} \}
$$
is chosen so that

a) the closed polydisc

$$
\overline{U_{\rho}}=\{z\in \mathbb{C}^n:|z_1|\leq e^{{\rho}_1},|z_2|\leq e^{{\rho}_2},|z_3|\leq e^{{\rho}_3} \}
$$

doesn't contain poles of $F(z)$;

b) parameters $t=(t_1,t_2)$ satisfy the inequalities  $|t_i|<e^{\langle q_i,\rho\rangle},i=1,2$.

Now let's form the matrix
 $A=\begin{pmatrix}
1&1&0\\
1&2&0\\
1&2&1
\end{pmatrix}$,
then  $A^{-1}=\begin{pmatrix}
2 & -1&0\\
-1& 1&0\\
0 &-1&1
\end{pmatrix}$. Using substitution $z^{A^{-1}}=w$, we get $z_1=w_1^2w_2^{-1},z_2=w_1^{-1}w_2^{1}w_3^{-1},z_3=w_3$. The denominator of the function $ F $ in variables $w$ looks like $1+w_1^2w_2^{-1}+w_1^{-1}w_2^{1}w_3^{-1}+w_3+w_1^{-1}w_2^{1}$ and it's Newton polytope is shown in Figure 2.

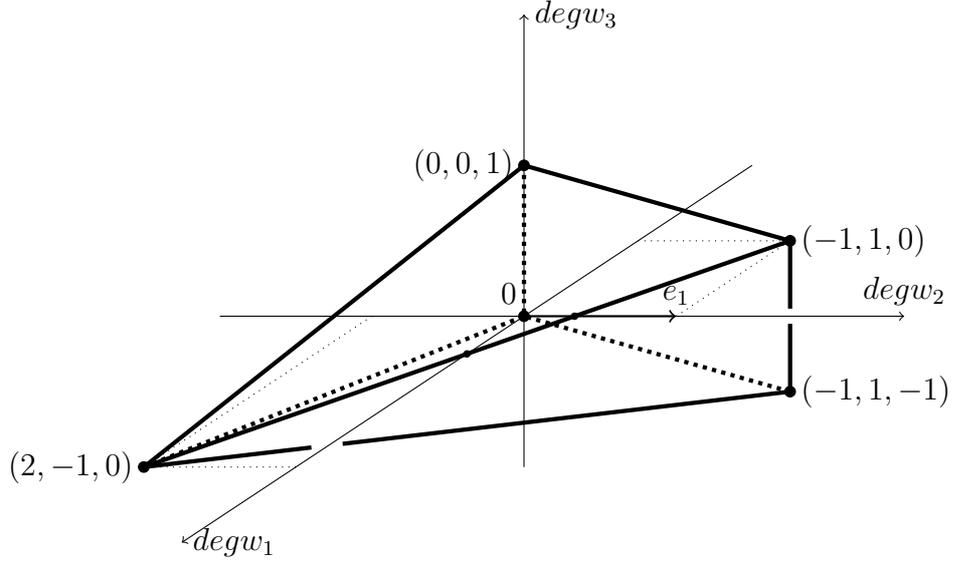
\begin{figure}[h!]
\begin{center}
\begin{tikzpicture}[domain=-4:4]
\draw [->](-4,0)--(5,0) node[above]{$deg w_2$};
\draw [->](3,2)--(-4.5,-3) node[right]{$deg w_1$};
\draw [->](0,-2)--(0,4) node[right]{$deg w_3$};
\draw [ultra thick](-5,-2)--(0,2);
\draw [ultra thick](0,2)--(3.5,1) ;
\draw [ultra thick](3.5,1)--(3.5,0.1) ;
\draw [ultra thick](3.5,-0.1)--(3.5,-1) ;
\draw [ultra thick](3.5-8.5/1.35,-1-1/1.35)--(-5,-2);
\draw [ultra thick](3.5,-1)--(-5+8.5/3.25,-2+1/3.25);
\draw [ultra thick](-5,-2)--(3.5,1) ;
\draw [dotted][ultra thick](-5,-2)--(0,0) ;
\draw [dotted][ultra thick](3.5,-1)--(0,0) ;
\draw [dotted][ultra thick](0,2)--(0,0) ;
\draw [dotted](1.5,1)--(3.5,1) ;
\draw [dotted](2,0)--(3.5,1) ;
\draw [dotted](-2,0)--(-5,-2) ;
\draw [dotted](-3,-2)--(-5,-2) ;
\draw [fill] (0,2) circle (2pt)node[left]{$(0,0,1)$} (-5,-2) circle (2pt)node[left]{$(2,-1,0)$} (3.5,1) circle (2pt)node[right]{$(-1,1,0)$} (3.5,-1) circle (2pt)node[right]{$(-1,1,-1)$} (0,0) circle (2pt);
\draw (-0.2,0.3) node{$ 0$};
\draw [fill] ( 0.666,0) circle (1.2pt) (-0.75,-0.5) circle (1.2pt);
\draw [thick,->](0,0)--(2,0) node[above]{$ e_1$};
\end{tikzpicture}
\caption{The Newton polytope of $Q(w^A)=1+w_1^2w_2^{-1}+w_1^{-1}w_2^{1}w_3^{-1}+w_3+w_1^{-1}w_2^{1}$.}
\end{center}
\end{figure} 

The integral (\ref{b8}) after the change of variables $z^{A^{-1}}=w$  has the following form:

\begin{equation*}
d_q(t_1,t_2)=\frac{1}{(2\pi i)^3}\int_{\Gamma_{\rho}}\frac{1}{1+w_1^2w_2^{-1}+w_1^{-1}w_2^{1}w_3^{-1}+w_3+w_1^{-1}w_2^{1}}\cdot
\frac{w_1\cdot w_2}{(w_1-t_1)(w_2-t_2)}\cdot
\frac{dw_1}{w_1}\frac{dw_2}{w_2}\frac{dw_3}{w_3}.
\end{equation*}

Applying the Cauchy formula with the variable $w_1$, we obtain the following representation of the diagonal

\begin{equation*}
d_q(t_1,t_2)=\frac{1}{(2\pi i)^2}\int_{\Gamma_{\rho'}}\frac{1}{1+t_1^2w_2^{-1}+t_1^{-1}w_2^{1}w_3^{-1}+w_3+t_1^{-1}w_2^{1}}\cdot
\frac{w_2}{w_2-t_2}\cdot
\frac{dw_2}{w_2}\frac{dw_3}{w_3}.
\end{equation*}
The Newton polytope of the denominator of the function $F(t_1,w_2,w_3)$ is shown in Figure~3.
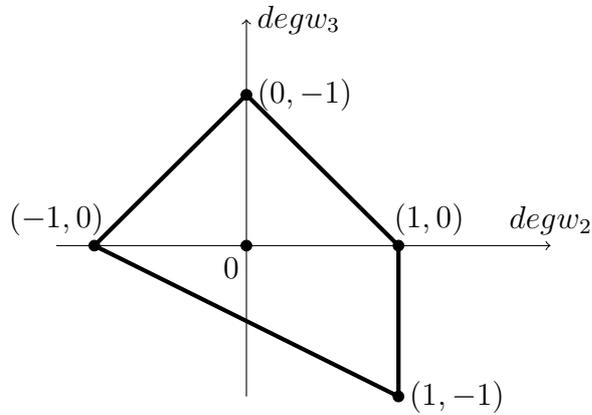
\begin{figure}[h]

\begin{center}
\begin{tikzpicture}[domain=-4:4]

\draw [->](5.5,0)--(12,0) node[above]{$deg w_2$};
\draw [->](8,-2)--(8,3) node[right]{$deg w_3$};
\draw (5.5,0.35)node{$(-1,0)$};
\draw (10.4,0.35)node{$(1,0)$};
\draw [fill] (8,0) circle (2pt)(6,0) circle (2pt)(8,2) circle (2pt)node[right]{$(0,-1)$}(10,-2) circle (2pt)node[right]{$(1,-1)$} (10,0) circle (2pt);
\draw [ultra thick](6,0)--(8,2) ;
\draw (7.8,-0.3) node{$ 0$};
\draw [ultra thick](8,2)--(10,0);
\draw [ultra thick](10,0)--(10,-2);
\draw [ultra thick](10,-2)--(6,0);
\end{tikzpicture}
\caption{The Newton polytope of $1+t_1^2w_2^{-1}+t_1^{-1}w_2^{1}w_3^{-1}+w_3+t_1^{-1}w_2^{1}$.}
\end{center}
\end{figure}

Integrating by the variable $ w_2 $, we obtain the representation
\begin{equation}\label{b9}
d_q(t_1,t_2)=\frac{1}{(2\pi i)}\int_{\Gamma_{\rho''}}\frac{1}{1+t_1^2t_2^{-1}+t_1^{-1}t_2^{1}w_3^{-1}+w_3+t_1^{-1}t_2^{1}}\cdot\frac{dw_3}{w_3}.
\end{equation}


Therefore the integral (\ref{b8}) admits the reduction to the one-dimensional integral (\ref{b9}) with the rational integrand. It is known (see Section~10.2 in \cite{Tsikh},) that such integral is an algebraic function in variables $t_1, t_2$.


\addcontentsline{toc}{section}{References}
	\begin{center}
	\renewcommand{\refname}
	\end{center}

\bigskip

\noindent {\bf Authors' addresses}: \noindent Artem Senashov, asenashov@mail.ru

\end{document}